\documentclass[10pt,twocolumn,twoside]{IEEEtran}
\IEEEoverridecommandlockouts

\usepackage{cite}
\usepackage{amsmath,amssymb,amsfonts}
\usepackage{graphicx}
\usepackage{textcomp}
\usepackage{balance}
\usepackage{xcolor}
\usepackage[utf8]{inputenc}
\usepackage{url}
\usepackage{eurosym}
\usepackage{mathtools}
\usepackage{algorithm, algpseudocode}
\usepackage{multirow}

\def\BibTeX{{\rm B\kern-.05em{\sc i\kern-.025em b}\kern-.08em
    T\kern-.1667em\lower.7ex\hbox{E}\kern-.125emX}}
\usepackage{lipsum}

\title{\LARGE \bf 
Unlocking the Informational Value of Marginal Costs for Exact Time Series Aggregation in Generation Expansion Planning
}

\author{Luca Santosuosso, Sonja Wogrin
\vspace{-0.5cm}
\thanks{ 
Luca Santosuosso and Sonja Wogrin are with the Institute of Electricity Economics and Energy Innovation, Graz University of Technology, and the Research Center ENERGETIC, 8010 Graz, Austria. Emails: {\tt luca.santosuosso@tugraz.at; wogrin@tugraz.at}.}
}

\begin{document}
\begingroup
\allowdisplaybreaks

\maketitle

\begin{abstract}
This paper addresses the generation expansion planning (GEP) problem, formulated as a mixed-integer linear programming model with intertemporal storage constraints. 
Being generally NP-hard, the problem’s computational complexity grows sharply with the planning horizon and the number of binary variables. 
While previous research has tackled this challenge using heuristic time series aggregation (TSA) methods, 
we propose a theoretically grounded marginal-cost-based TSA, designed to construct an aggregated model that preserves the active constraints of its full-scale counterpart, thereby explicitly targeting exact temporal aggregation. 
This TSA method is embedded within solution algorithms that iteratively refine theoretically validated bounds on the maximum error introduced by the temporal aggregation relative to conventional full-scale optimization, thus offering a formal performance guarantee to the decision-maker. 
Numerical results highlight the computational advantages of the proposed algorithms, which notably recover tractability whereas full-scale optimization proves intractable.
\end{abstract}

\begin{IEEEkeywords}
Exact time series aggregation, generation expansion planning, energy storage, mixed-integer programming
\end{IEEEkeywords}

\section{Introduction}
\label{sec:introduction}

\subsection{Literature Review}

Generation expansion planning (GEP) is a cornerstone problem in energy system optimization,
aiming to determine the optimal mix and sizing of power generation units required to meet future electricity demand
while minimizing both capital investment and operational costs \cite{BARINGO2023101105}.
Recently, the large-scale integration of variable renewable energy sources (vRES),
the growing deployment of flexible technologies (e.g., energy storage systems), and the evolution of regulatory and market designs
have significantly increased the complexity of GEP models \cite{koltsaklis2018state}.
Modern GEP formulations must therefore capture nonlinear system dynamics, alongside intricate intertemporal coupling constraints, over long-term planning horizons \cite{motta2024survey}.
This has driven a steady evolution of GEP models from basic linear programming (LP) formulations \cite{hobbs1995optimization} toward increasingly complex nonconvex formulations \cite{franco2014mixed},
most notably mixed-integer linear programming (MILP) models \cite{lara2018deterministic}.
Critically, MILP formulations of the GEP problem are strongly NP-hard in general \cite{GODERBAUER2019343}, making the trade-off between modeling fidelity and computational tractability a central challenge \cite{li2022mixed}.

To address this trade-off, time series aggregation (TSA) methods are commonly employed \cite{8369128}.
By leveraging clustering techniques, these methods condense the input time series of the GEP problem,
such as energy demand profiles or vRES generation forecasts,
into a temporally aggregated counterpart of the original full-scale GEP model,
defined over a reduced set of representative time periods (or clusters) \cite{granell2014impacts}.

Traditional \textit{a priori} TSA methods rely on standard clustering techniques,
such as k-means \cite{sarajpoor2023time}, k-medoids \cite{SCHUTZ2018570}, and hierarchical clustering \cite{8017598},
to select representative periods based exclusively on the statistical properties of the input time series.
Despite their widespread adoption in GEP \cite{munoz2016new}, the accuracy of these methods is highly case-dependent and must therefore be validated individually for each application \cite{KOTZUR2018474}.
This limitation has motivated growing interest in establishing performance guarantees for aggregated models obtained via TSA \cite{11305411},
typically formulated as bounds on the approximation error relative to the original full-scale model \cite{rogers1991aggregation}.
Although such bounds are valuable for enabling decision-makers to assess the accuracy of an aggregated model \cite{tso2020hierarchical},
they do not allow a priori TSA-based models to systematically attain a prescribed level of accuracy.
This shortcoming reflects a fundamental limitation of a priori methods:
they implicitly assume that an accurate statistical representation of the full-scale model’s input space yields a correspondingly accurate approximation of its output solution space.
In practice, however, this assumption often fails, as the model outputs depend not only on the input time series but also critically on the underlying mathematical structure of the full-scale model \cite{sun2019data}.

The recognition of this limitation has led to the recent emergence of the \textit{a posteriori} TSA paradigm \cite{hoffmann2020review}.
The key distinction of this paradigm from conventional a priori TSA is its use of clustering features derived directly from the full-scale model,
rather than relying solely on statistical features extracted from the input time series \cite{li2022representative}.
This paradigm has proven highly effective for GEP \cite{poncelet2016selecting}, particularly in the presence of intertemporal coupling constraints (e.g., due to storage or ramping) \cite{tejada2018enhanced},
which significantly complicate TSA by requiring preservation of temporal chronology to maintain consistency with the full-scale model dynamics \cite{kotzur2018time}.
Existing studies have explored a posteriori methods based on clustering features derived from approximations of the full-scale primal solution \cite{zhang2022model} or the dual solution \cite{cardona2024enhancing}.
Nevertheless, a posteriori TSA remains underdeveloped, as the current research has yet to conclusively determine which clustering features most effectively lead to accurate aggregated models.

In this regard, the theoretical result of \cite{10037240} offers a strong conceptual foundation by demonstrating that,
if the decision-maker had prior knowledge of the full-scale model’s active constraints,
TSA could be performed \textbf{exactly},
with the aggregated formulation reproducing both the optimal solution space and the optimal objective function value of its full-scale counterpart.
Notably, \cite{klatzer2025towards} further extends this result to optimization models featuring intertemporal storage constraints.

While perfect knowledge of the full-scale model’s active constraint sets during TSA is clearly unattainable, the theoretical result of \cite{10037240} provides a valuable insight:
any clustering feature that reliably captures constraint activation in the full-scale model can be exploited within an a posteriori TSA method to achieve exact temporal aggregation.
This highlights a significant, yet largely untapped, opportunity to dramatically reduce the computational complexity of nonconvex GEP models, without compromising solution accuracy.

\subsection{Research Gaps and Contributions}
The literature review reveals the following research gaps:
\begin{itemize}
    \item Traditional \emph{a priori} TSA methods rely on the frequently violated assumption that an accurate representation of the full-scale model’s input space is sufficient to accurately capture its solution space \cite{sun2019data}.
    In contrast, \emph{a posteriori} TSA methods offer greater potential for constructing high-fidelity aggregated models by exploiting information extracted directly from the full-scale model to enrich the clustering feature space \cite{hoffmann2020review}. 
    Although prior studies have proposed deriving clustering features from approximations of the full-scale primal \cite{zhang2022model} and dual \cite{cardona2024enhancing} solutions,
    it remains unclear which full-scale model information should be leveraged, and how clustering features should be systematically constructed, to improve the accuracy of temporally aggregated models.
    This challenge is particularly pronounced for MILP optimization models, where dual information is generally ill-defined.
    \item Recent analyses \cite{10037240,klatzer2025towards} demonstrate that exact temporal aggregation can be achieved when the active constraints of the full-scale model are preserved in the aggregated model.
    These results provide clear theoretical guidance on which aspects of the full-scale model are most critical for ensuring aggregation accuracy and, therefore, should be captured by TSA methods.
    However, because the set of active constraints is unknown during the TSA process, it remains unclear how this theoretical insight can be translated into practical methods that effectively reduce GEP models while maintaining high solution accuracy and formal performance guarantees.
\end{itemize}

This paper seeks to address these research gaps through the following \textbf{key contributions}:
\begin{itemize}
    \item We formulate the full-scale GEP model as a MILP problem encompassing multiple generation technologies, including both vRES and thermal power, and subject to intertemporal storage constraints.
    We then construct the corresponding temporally aggregated model and show that the full-scale model’s marginal costs serve as a reliable proxy for identifying the sets of active constraints in the full-scale model, enabling their use as clustering features to achieve exact temporal aggregation.
    \item Acknowledging that system marginal costs are not directly available during temporal aggregation,
    we introduce a practical a posteriori TSA method that leverages marginal cost estimates obtained from the dual space of the full-scale model’s linear relaxation
    to guide the construction of the temporally aggregated model.
    \item We develop two iterative solution algorithms that systematically reduce the temporal dimensionality of the full-scale MILP GEP model while providing a formal performance guarantee
    in the form of theoretically validated bounds on the approximation error introduced by temporal aggregation.
    Notably, both algorithms deliver feasible solutions for the full-scale model at each iteration,
    and the aggregated model is progressively refined toward exact temporal aggregation through the proposed marginal-costs-based a posteriori TSA method.
\end{itemize}

Finally, the proposed algorithms are evaluated using real-world data from the ENTSO-E Transparency Platform \cite{hirth2018entso}.

The remainder of the paper is organized as follows:
Section~\ref{sec:methodology} details the proposed methodology,
Section~\ref{sec:results} presents the numerical results,
and Section~\ref{sec:conclusion} concludes the study.

\section{Methodology}
\label{sec:methodology}
This section presents the proposed methodology.
Subsection~\ref{subsec:full_model} introduces the full-scale GEP model,
while Subsection~\ref{subsec:aggregated_model} presents its temporally aggregated counterpart.
Subsection~\ref{subsec:marginal_cost_clustering} details the proposed marginal-cost-based a posteriori TSA method,
which is employed in the algorithms of Subsection~\ref{subsec:algorithms} to iteratively refine bounds on the optimal objective function value of the original full-scale model.

In the following, sets and vectors are denoted in boldface (e.g., $\boldsymbol{z}$). The notation $|\cdot|$ denotes the cardinality of a set.

\subsection{The Full-Scale Model}
\label{subsec:full_model}
The goal is to determine the optimal mix and sizing of generation units that minimize capital investment and operational costs while meeting the energy demand.
Specifically, we consider an asset portfolio comprising vRES, namely wind and solar power plants, thermal units, and energy storage systems.
In the following formulation, we omit network constraints, consistent with the operation of generation portfolios that do not participate in the grid management, such as commercial virtual power plants or generation companies \cite{BARINGO2023101105}.
Future work will extend the formulation to incorporate these constraints.

Let $\boldsymbol{G}$, $\boldsymbol{N}$, and $\boldsymbol{T}$ denote the sets of generators (indexed by $g$), energy storage systems (indexed by $n$), and time steps (indexed by $t$), respectively.
The operational costs of generator $g$ and non-supplied energy are denoted by
$C^{\mathrm{pg}}_g$ and $C^{\mathrm{ns}}$ (\euro/MWh), respectively.
For each storage unit $n$,
$C^{\mathrm{pc}}_n$ and $C^{\mathrm{pd}}_n$ denote the charging and discharging costs (\euro/MWh).
The investment costs (\euro/MW) are 
$C^{\mathrm{g,inv}}_g$ for generators and $C^{\mathrm{s,inv}}_n$ for storage units.
The input time series comprise the generator capacity factors $F_{g,t}$ and the energy demand $D_t$ (MWh).

We denote by
$x^{\mathrm{g}}_g$ and $x^{\mathrm{s}}_n$
the installed capacities (MW) of generator $g$ and storage unit $n$, respectively.
The binary variables $b^{\mathrm{g}}_g$ and $b^{\mathrm{s}}_n$
enforce that installed capacities are either zero or lie within
$\left[\underline{X}^{\mathrm{g}}_g,\overline{X}^{\mathrm{g}}_g\right]$
and $\left[\underline{X}^{\mathrm{s}}_n,\overline{X}^{\mathrm{s}}_n\right]$, respectively.

The operational variables include the power output $p^{\mathrm{g}}_{g,t}$ (MW) of generator $g$ at time $t$,
the non-supplied energy $e^{\mathrm{ns}}_{t}$ (MWh), 
the charging and discharging power of storage unit $n$, denoted by $p^{\mathrm{c}}_{n,t}$ and $p^{\mathrm{d}}_{n,t}$ (MW), respectively,
and its stored energy $e^{\mathrm{s}}_{n,t}$ (MWh) at time $t$.
The charging and discharging efficiencies of storage unit $n$ are denoted by $\eta^{\mathrm{c}}_n$ and $\eta^{\mathrm{d}}_n$, respectively,
$\rho_{n}$ denotes its energy-to-power ratio (h),
while $E^0_n$ denotes its initial stored energy (MWh).

We group the decision variables of the full-scale model in the set $\boldsymbol{z}$, defined as:
\begin{equation*}
    \boldsymbol{z} \!\coloneqq\! \left\{ x^{\mathrm{g}}_g, x^{\mathrm{s}}_n, b^{\mathrm{g}}_g, b^{\mathrm{s}}_n, p^{\mathrm{g}}_{g,t}, e^{\mathrm{s}}_{n,t}, p^{\mathrm{c}}_{n,t}, p^{\mathrm{d}}_{n,t}, e^\mathrm{ns}_{t} \right\}_{g \in \boldsymbol{G}, n \in \boldsymbol{N}, t \in \boldsymbol{T}}.
\end{equation*}

The objective function of the GEP problem is defined as
\begin{align}
\label{full_model:obj}
J(\boldsymbol{z}) \! \coloneqq \! & \sum_{g \in \boldsymbol{G}} C^\mathrm{g,inv}_g x^{\mathrm{g}}_g + \sum_{n \in \boldsymbol{N}} C^\mathrm{s,inv}_n x^{\mathrm{s}}_n + \sum_{t \in \boldsymbol{T}} \sum_{g \in \boldsymbol{G}} C^\mathrm{pg}_g p^\mathrm{g}_{g,t} \Delta \nonumber\\
& + \sum_{t \in \boldsymbol{T}} \! \left(\sum_{n \in \boldsymbol{N}} \!\! \left(C^{\mathrm{pc}}_n p^{\mathrm{c}}_{n,t} + C^{\mathrm{pd}}_n p^{\mathrm{d}}_{n,t}\right) \Delta + C^\mathrm{ns} e^\mathrm{ns}_t \!\right),
\end{align}
which represents the sum of both capital investment and operational costs over the planning horizon $\boldsymbol{T}$.

The \textbf{full-scale GEP model}, formulated as a MILP problem, is defined over $\boldsymbol{T}$ with sampling time $\Delta$ (h) as follows:
\begin{subequations}
\label{full_model}
\begin{align}
\min_{\boldsymbol{z}} \quad & J\left(\boldsymbol{z}\right)\\
\textrm{s.t.} \quad & \sum_{g \in \boldsymbol{G}} p^{\mathrm{g}}_{g, t} \Delta + \!\! \sum_{n \in \boldsymbol{N}} \! \left(p^{\mathrm{d}}_{n,t} - p^{\mathrm{c}}_{n,t}\right) \! \Delta + e^\mathrm{ns}_t \!=\! D_t, \; \forall t, \label{full_model:power_bal}\\
& e^{\mathrm{s}}_{n,t+1} = e^{\mathrm{s}}_{n,t} + \left(\eta^{\mathrm{c}}_n p^{\mathrm{c}}_{n,t} - \eta^{\mathrm{d}}_n p^{\mathrm{d}}_{n,t}\right) \Delta, \nonumber\\
& \qquad\qquad\qquad\qquad\quad\;\;\; \forall n, \forall t \in \boldsymbol{T} \setminus \left\{|\boldsymbol{T}|-1\right\}, \label{full_model:sto_dyn}\\
& e^{\mathrm{s}}_{n,0} = E^0_n, \; \forall n, \label{full_model:sto_init}\\
& 0 \leq p^{\mathrm{g}}_{g,t} \leq F_{g,t} \, x^{\mathrm{g}}_g, \; \forall g, \forall t, \label{full_model:gen_lim}\\
& 0 \leq e^{\mathrm{s}}_{n,t} \leq x^{\mathrm{s}}_n \, \rho_n, \; \forall n, \forall t, \label{full_model:sto_lim}\\
& 0 \leq p^{\mathrm{c}}_{n,t} \leq x^{\mathrm{s}}_n, \; \forall n, \forall t, \label{full_model:sto_c_lim}\\
& 0 \leq p^{\mathrm{d}}_{n,t} \leq x^{\mathrm{s}}_n, \; \forall n, \forall t, \label{full_model:sto_d_lim}\\
& b^{\mathrm{s}}_n \, \underline{X}^{\mathrm{s}}_n \leq x^{\mathrm{s}}_{n} \leq b^{\mathrm{s}}_n \, \overline{X}^{\mathrm{s}}_n, \; \forall n, \label{full_model:sto_inv_lim}\\
& b^{\mathrm{g}}_g \, \underline{X}^{\mathrm{g}}_g \leq x^{\mathrm{g}}_{g} \leq b^{\mathrm{g}}_g \, \overline{X}^{\mathrm{g}}_g, \; \forall g, \label{full_model:gen_inv_lim}\\
& e^\mathrm{ns}_t \geq 0, \; \forall t, \label{full_model:nse}\\
& b^{\mathrm{g}}_g \in \left\{0,1\right\}, \; \forall g, \label{full_model:bin_g}\\
& b^{\mathrm{s}}_n \in \left\{0,1\right\}, \; \forall n. \label{full_model:bin_n}
\end{align}
\end{subequations}
In \eqref{full_model},
the constraints \eqref{full_model:power_bal} enforce the energy balance;
\eqref{full_model:sto_dyn} and \eqref{full_model:sto_init} characterize the storage dynamics;
\eqref{full_model:sto_inv_lim} and \eqref{full_model:gen_inv_lim} impose capacity investment limits;
and \eqref{full_model:gen_lim}--\eqref{full_model:sto_d_lim} together with \eqref{full_model:nse}--\eqref{full_model:bin_n} define the bounds on the decision variables.

\subsection{The Temporally Aggregated Model}
\label{subsec:aggregated_model}

Mixed-integer GEP problems are known to be strongly NP-hard \cite{GODERBAUER2019343}.
As the problem dimension increases with the number of binary variables and time steps,
\eqref{full_model} may become computationally intractable.
To mitigate this, TSA can be employed to construct an aggregated counterpart of \eqref{full_model},
defined over a reduced set of representative periods, denoted by $\boldsymbol{K}$ and indexed by $k$.
When $|\boldsymbol{K}| \ll |\boldsymbol{T}|$, the aggregated model offers a significant computational advantage over the full-scale model.

Let $\boldsymbol{T}_k$ denote the set of consecutive time steps $t \in \boldsymbol{T}$ assigned to the $k$-th cluster via TSA.
We define $K := |\boldsymbol{K}|$ and $T_k := |\boldsymbol{T}_k|$.
The aggregated counterparts of $F_{g,t}$ and $D_t$ are computed as

\vspace{-0.1cm}
\begin{minipage}{0.50\columnwidth}
\begin{equation}
\label{prop:main_result_eq1}
\hat{F}_{g,k} \!\coloneqq\! \sum_{t \in \boldsymbol{T}_k} \! \frac{F_{g,t}}{T_k}, \; \forall g, \forall k,
\end{equation}
\end{minipage}%
\hfill
\begin{minipage}{0.45\columnwidth}
\begin{equation}
\label{prop:main_result_eq2}
\hat{D}_{k} \coloneqq \sum_{t \in \boldsymbol{T}_k} \frac{D_{t}}{T_k}, \; \forall k.
\end{equation}
\end{minipage}
\vspace{0.01cm}

We group the variables of the aggregated model as follows:
\begin{equation*}
    \boldsymbol{\hat{z}} \!\coloneqq\! \left\{ \hat{x}^{\mathrm{g}}_g, \hat{x}^{\mathrm{s}}_n, \hat{b}^{\mathrm{g}}_g, \hat{b}^{\mathrm{s}}_n, \hat{p}^{\mathrm{g}}_{g,k}, \hat{e}^{\mathrm{s}}_{n,k}, \hat{p}^{\mathrm{c}}_{n,k}, \hat{p}^{\mathrm{d}}_{n,k}, \hat{e}^\mathrm{ns}_{k} \right\}_{g \in \boldsymbol{G}, n \in \boldsymbol{N}, k \in \boldsymbol{K}}.
\end{equation*}

The aggregated counterpart of $J(\boldsymbol{z})$ in \eqref{full_model:obj} is
\begin{align}
\label{agg_model:obj}
\hat{J}(\boldsymbol{\hat{z}}) \! \coloneqq \! & \sum_{g \in \boldsymbol{G}} \! C^\mathrm{g,inv}_g \hat{x}^{\mathrm{g}}_g + \!\! \sum_{n \in \boldsymbol{N}} \! C^\mathrm{s,inv}_n \hat{x}^{\mathrm{s}}_n \!+ \!\! \sum_{k \in \boldsymbol{K}} \sum_{g \in \boldsymbol{G}} \! C^\mathrm{pg}_g \hat{p}^\mathrm{g}_{g,k} T_k \Delta \nonumber\\
& + \!\! \sum_{k \in \boldsymbol{K}} \!\! \left(\sum_{n \in \boldsymbol{N}} \!\! \left(C^{\mathrm{pc}}_n \hat{p}^{\mathrm{c}}_{n,k} \!+\! C^{\mathrm{pd}}_n \hat{p}^{\mathrm{d}}_{n,k}\right) \!\Delta \!+\! C^\mathrm{ns} \hat{e}^\mathrm{ns}_{k} \!\!\right) \! T_k.
\end{align}

The \textbf{temporally aggregated GEP model}, formulated as a MILP problem, is defined over $\boldsymbol{K}$ as follows:
\begin{subequations}
\label{agg_model}
\begin{align}
\min_{\boldsymbol{\hat{z}}} \quad & \hat{J}\left(\boldsymbol{\hat{z}}\right)\\
\textrm{s.t.} \quad & \!\sum_{g \in \boldsymbol{G}} \hat{p}^{\mathrm{g}}_{g,k} \Delta \!+ \!\! \sum_{n \in \boldsymbol{N}} \!\! \left(\hat{p}^{\mathrm{d}}_{n,k} - \hat{p}^{\mathrm{c}}_{n,k}\right) \! \Delta \!+\! \hat{e}^\mathrm{ns}_{k} \!=\! \hat{D}_{k}, \; \forall k, \label{agg_model:power_bal}\\
& \hat{e}^{\mathrm{s}}_{n,k+1} = \hat{e}^{\mathrm{s}}_{n,k} + \left(\eta^{\mathrm{c}}_n \hat{p}^{\mathrm{c}}_{n,k} - \eta^{\mathrm{d}}_n \hat{p}^{\mathrm{d}}_{n,k}\right) T_{k} \, \Delta, \nonumber\\
& \qquad\qquad\qquad\qquad\quad\;\;\; \forall n, \forall k \in \boldsymbol{K} \setminus \left\{K-1\right\}, \label{agg_model:sto_dyn}\\
& \hat{e}^{\mathrm{s}}_{n,0} = E^0_n, \; \forall n, \label{agg_model:sto_init}\\
& 0 \leq \hat{p}^{\mathrm{g}}_{g,k} \leq \hat{F}_{g,k} \, \hat{x}^{\mathrm{g}}_g, \; \forall g, \forall k, \label{agg_model:gen_lim}\\
& 0 \leq \hat{e}^{\mathrm{s}}_{n,k} \leq \hat{x}^{\mathrm{s}}_n \, \rho_n, \; \forall n, \forall k, \label{agg_model:sto_lim}\\
& 0 \leq \hat{p}^{\mathrm{c}}_{n,k} \leq \hat{x}^{\mathrm{s}}_n, \; \forall n, \forall k, \label{agg_model:sto_c_lim}\\
& 0 \leq \hat{p}^{\mathrm{d}}_{n,k} \leq \hat{x}^{\mathrm{s}}_n, \; \forall n, \forall k, \label{agg_model:sto_d_lim}\\
& \hat{b}^{\mathrm{s}}_n \, \underline{X}^{\mathrm{s}}_n \leq \hat{x}^{\mathrm{s}}_{n} \leq \hat{b}^{\mathrm{s}}_n \, \overline{X}^{\mathrm{s}}_n, \; \forall n, \label{agg_model:sto_inv_lim}\\
& \hat{b}^{\mathrm{g}}_g \, \underline{X}^{\mathrm{g}}_g \leq \hat{x}^{\mathrm{g}}_{g} \leq \hat{b}^{\mathrm{g}}_g \, \overline{X}^{\mathrm{g}}_g, \; \forall g, \label{agg_model:gen_inv_lim}\\
& \hat{e}^\mathrm{ns}_k \geq 0, \; \forall k, \label{agg_model:nse}\\
& \hat{b}^{\mathrm{g}}_g \in \left\{0,1\right\}, \; \forall g, \label{agg_model:bin_g}\\
& \hat{b}^{\mathrm{s}}_n \in \left\{0,1\right\}, \; \forall n. \label{agg_model:bin_n}
\end{align}
\end{subequations}
In \eqref{agg_model},
the constraints \eqref{agg_model:power_bal}–\eqref{agg_model:bin_n}
represent the temporally aggregated counterparts of 
\eqref{full_model:power_bal}–\eqref{full_model:bin_n}, respectively.

\subsection{Clustering Based on Marginal Costs}
\label{subsec:marginal_cost_clustering}


The construction of the aggregated model \eqref{agg_model} requires clustering to identify representative periods within the planning horizon $\boldsymbol{T}$.
The intertemporal constraints in \eqref{full_model} render standard methods (e.g., k-means) unsuitable, as they neglect temporal chronology; we therefore use a \emph{sliding-window clustering} to group consecutive time steps by clustering feature similarity.

Let $\{\boldsymbol{a}_t\}_{t \in \boldsymbol{T}}$ be the set of clustering features.
We evaluate the similarity between consecutive elements of $\{\boldsymbol{a}_t\}$,
assigning $t$ to cluster $k$ if
\begin{equation}
\label{TSACondition}
\|\boldsymbol{a}_t - \boldsymbol{\mu}_k\|_2 \leq \zeta,
\end{equation}
where $\boldsymbol{\mu}_k$ denotes the centroid of cluster $k$, defined as $\boldsymbol{\mu}_k \coloneqq \frac{1}{T_k} \sum_{t \in \boldsymbol{T}_k} \boldsymbol{a}_t$, and $\zeta$ is a user-defined similarity threshold.
If \eqref{TSACondition} is not met for $\boldsymbol{a}_t$, a new cluster is initialized starting at $t$.

The selection of clustering features is a key determinant of the quality of the aggregated model \eqref{agg_model}.
Conventional a priori TSA methods construct these features exclusively from the input time series of the GEP problem.
Although this strategy preserves the problem input space,
it offers no guarantee that the aggregated model will accurately reproduce the outputs
(e.g., the optimal objective function value or the optimal solution space)
of the full-scale model \cite{sun2019data}.
In general, a decision-maker is primarily concerned with preserving accuracy in the GEP model’s output space rather than its input space.

To address this issue, we build on the theoretical results of \cite{10037240} and \cite{klatzer2025towards},
which establish that an aggregated model reproduces the exact outputs of its full-scale counterpart when it captures the active constraints of the full-scale model at optimality.
In the full-scale model \eqref{full_model},
the constraints \eqref{full_model:gen_lim}--\eqref{full_model:sto_d_lim} become active when a generator reaches its operational limits,
as dictated by \eqref{full_model:sto_inv_lim} and \eqref{full_model:gen_inv_lim},
thereby triggering a transition in the marginal generator that sets the system’s marginal cost,
which represents the sensitivity of the optimal total system cost to a marginal change in the energy demand.
Hence, marginal costs provide an effective proxy for constraint activation and can be used as clustering features to guide the construction of \eqref{agg_model} toward exact temporal aggregation.

\subsection{Practical Solution Algorithms}
\label{subsec:algorithms}

\begin{algorithm}[t]
\caption{Marginal Cost–Based Time Series Aggregation with Bounded Objective Function Error}\label{alg1}
\begin{algorithmic}[1]
\Require Parameters of \eqref{full_model}, optimality threshold $\epsilon^\mathrm{thr}$, maximum number of iterations $I \coloneqq |\boldsymbol{I}|$, and similarity threshold $\zeta$.

\Ensure Objective function bounds ${{J^\mathrm{UB}}}^\star$ and ${{J^\mathrm{LB}}}^\star$.

\Statex $i \gets 1$; $\epsilon^i \gets +\infty$; $J^{\mathrm{LB}^{i}} \gets -\infty$; $J^{\mathrm{UB}^{i}} \gets +\infty$;

\While{$\epsilon^i > \epsilon^{\mathrm{thr}}$ and $i \leq I$}

\State \parbox[t]{\dimexpr\linewidth-\algorithmicindent}{Estimate marginal costs according to Steps I and II of Subsection~\ref{subsec:algorithms};}

\State \parbox[t]{\dimexpr\linewidth-\algorithmicindent}{$\left\{\boldsymbol{T}_k^i\right\}_{k \in \boldsymbol{K}^i} \gets$ Sliding-window clustering of Subsection~\ref{subsec:marginal_cost_clustering} using threshold $\zeta$ and estimated marginal costs as features;}

\State $\boldsymbol{\hat{z}}^\star \gets$ Solve the aggregated model \eqref{agg_model} for $\left\{\boldsymbol{T}_k^i\right\}_{k \in \boldsymbol{K}^i}$;

\State $J^{\mathrm{LB}^{i+1}} \gets \mathrm{max} \left(\hat{J}\left(\boldsymbol{\hat{z}}^\star\right), J^{\mathrm{LB}^{i}}\right)$;

\vspace{1pt}

\State \parbox[t]{\dimexpr\linewidth-\algorithmicindent}{$\tilde{J} \gets$ Solve the full-scale model \eqref{full_model} with investment variable values fixed to those in $\boldsymbol{\hat{z}}^\star$ (LP model);}

\vspace{2pt}

\State $J^{\mathrm{UB}^{i+1}} \gets \mathrm{min} \left(\tilde{J}, J^{\mathrm{UB}^{i}}\right)$;

\vspace{2pt}

\State $\epsilon^{i+1} \gets 100 \, \frac{{{J^\mathrm{UB}}^{i + 1}} - {{J^\mathrm{LB}}^{i + 1}}}{{{J^\mathrm{UB}}^{i + 1}}}$;

\State $i \gets i+1$;

\EndWhile

\State $J^{\mathrm{LB}^\star} \gets J^{\mathrm{LB}^i}$ and $J^{\mathrm{UB}^\star} \gets J^{\mathrm{UB}^i}$;

\end{algorithmic}
\end{algorithm}

While Subsection \ref{subsec:marginal_cost_clustering} highlighted the potential of marginal cost-based clustering for exact TSA, two fundamental challenges remain. First, the nonconvexity of the full-scale model~\eqref{full_model} precludes the computation of marginal costs via duality theory. Second, even if marginal costs were obtainable, their evaluation would require solving the full-scale model itself, thereby negating the computational benefits of TSA.

To overcome these, we propose the following marginal cost estimation procedure, executed over $i \in \boldsymbol{I}$ iterations.

\textbf{Step I}.
For each month in the GEP horizon, randomly sample $i$ days. Next, solve a surrogate linear model obtained by relaxing the binary variables in~\eqref{full_model} to continuous variables in the interval $[0,1]$. This surrogate model is defined at full-scale temporal resolution but only over the set of sampled days.
For instance, if $i=2$ and the full-scale model is formulated at hourly resolution, the surrogate model spans $576$ time steps (i.e., 2 sampled days per month over 12 months, with 24 hours per day).
This optimization step yields an estimate of the marginal cost trajectory for each sampled day.

\textbf{Step II}.
For each month, all non-sampled days are assigned, with equal probability, to one of the $i$ sampled days of that month. 
Each non-sampled day is assumed to share the same marginal cost trajectory as the sampled day to which it is assigned. 
For example, if a particular day is assigned to a sampled day $d$, the hourly marginal costs of $d$ are used as a proxy for the marginal costs of that day.
This step yields a marginal cost estimate over the full-scale horizon $\boldsymbol{T}$, which is then used as a feature in the sliding-window clustering of Subsection~\ref{subsec:marginal_cost_clustering} to construct the aggregated model~\eqref{agg_model}.

\begin{algorithm}[t]
\caption{Adaptive Marginal Cost–Based Time Series Aggregation with Bounded Objective Function Error}\label{alg2}
\begin{algorithmic}[1]
\Require Parameters of \eqref{full_model}, optimality threshold $\epsilon^\mathrm{thr}$, maximum number of iterations $I \coloneqq |\boldsymbol{I}|$, and similarity threshold $\zeta$.

\Ensure Objective function bounds ${{J^\mathrm{UB}}}^\star$ and ${{J^\mathrm{LB}}}^\star$.

\State $i \gets 1$; $\epsilon^i \gets +\infty$; $J^{\mathrm{LB}^{i}} \gets -\infty$; $J^{\mathrm{UB}^{i}} \gets +\infty$;

\While{$\epsilon^i > \epsilon^{\mathrm{thr}}$ and $i \leq I$}

\If{$i = 1$}

\vspace{2pt}

\State \parbox[t]{\dimexpr\linewidth-\algorithmicindent}{$\boldsymbol{\hat{\lambda}}^{i} \gets$ Estimate marginal costs according to Steps \\ I and II of Subsection~\ref{subsec:algorithms};}

\Else

\State \parbox[t]{\dimexpr\linewidth-\algorithmicindent}{$\boldsymbol{\hat{\lambda}}^{i} \gets$ Estimate marginal costs based on the \\
discrepancies between $\boldsymbol{\tilde{\lambda}}^{i-1}$ and $\boldsymbol{\hat{\lambda}}^{i-1}$,\\
as detailed in Subsection~\ref{subsec:algorithms};}

\vspace{3pt}

\EndIf

\State \parbox[t]{\dimexpr\linewidth-\algorithmicindent}{$\left\{\boldsymbol{T}_k^i\right\}_{k \in \boldsymbol{K}^i} \gets$ Sliding-window clustering of Subsection~\ref{subsec:marginal_cost_clustering} using threshold $\zeta$ and $\boldsymbol{\hat{\lambda}}^{i}$ as features;}

\State $\boldsymbol{\hat{z}}^\star \gets$ Solve the aggregated model \eqref{agg_model} for $\left\{\boldsymbol{T}_k^i\right\}_{k \in \boldsymbol{K}^i}$;

\State $J^{\mathrm{LB}^{i+1}} \gets \mathrm{max} \left(\hat{J}\left(\boldsymbol{\hat{z}}^\star\right), J^{\mathrm{LB}^{i}}\right)$;

\vspace{1pt}

\State \parbox[t]{\dimexpr\linewidth-\algorithmicindent}{$\left\{ \boldsymbol{\tilde{\lambda}}^{i}, \tilde{J} \right\} \gets$ Solve the full-scale model \eqref{full_model} with investment variable values fixed to those in $\boldsymbol{\hat{z}}^\star$ (LP model);}

\vspace{2pt}

\State $J^{\mathrm{UB}^{i+1}} \gets \mathrm{min} \left(\tilde{J}, J^{\mathrm{UB}^{i}}\right)$;

\vspace{2pt}

\State $\epsilon^{i+1} \gets 100 \, \frac{{{J^\mathrm{UB}}^{i + 1}} - {{J^\mathrm{LB}}^{i + 1}}}{{{J^\mathrm{UB}}^{i + 1}}}$;

\State $i \gets i+1$;

\EndWhile

\State $J^{\mathrm{LB}^\star} \gets J^{\mathrm{LB}^i}$ and $J^{\mathrm{UB}^\star} \gets J^{\mathrm{UB}^i}$;

\end{algorithmic}
\end{algorithm}

Importantly, irrespective of the clustering technique employed, the aggregated model \eqref{agg_model} provides, by construction, a lower bound on the optimal objective function value of the full-scale model \eqref{full_model}~\cite{santosuosso2025we}.
Fixing the investment variables in \eqref{full_model} to those computed via \eqref{agg_model} reduces \eqref{full_model} to a linear program, whose solution yields an upper bound on the full-scale optimal objective.
At each iteration $i \in \boldsymbol{I}$, the relative difference between the upper and lower bounds, denoted $J^{\mathrm{UB}^{i}}$ and $J^{\mathrm{LB}^{i}}$, respectively, defines the optimality gap $\epsilon^i$, indicating that the MILP GEP solution derived via TSA is at most $\epsilon^i$ suboptimal.
Convergence is achieved when $\epsilon^i$ falls below a prescribed threshold $\epsilon^{\mathrm{thr}}$.
This iterative procedure, which combines marginal cost-based a posteriori TSA with bound-based error evaluation, constitutes the proposed Algorithm~\ref{alg1}.

Algorithm~\ref{alg1} refines the aggregated model by progressively increasing the number of sampled days used to estimate marginal costs.
However, the optimization steps employed to compute the objective function bounds also generate additional, underexploited information.
Specifically, solving~\eqref{full_model} with fixed investment variables to obtain an upper bound on the objective function produces short-run marginal cost estimates, which can be used to enhance the TSA as follows.

Let $\boldsymbol{\tilde{\lambda}}^i$ denote the marginal costs obtained from the optimization step yielding $J^{\mathrm{UB}^i}$ (i.e., the short-run marginal costs),
and let $\boldsymbol{\hat{\lambda}}^i$ denote those used to construct the aggregated model at iteration $i$ (i.e., the long-run marginal costs).
Algorithm~\ref{alg2} introduces an adaptive alternative to Algorithm~\ref{alg1}
executing the above Steps I and II at iteration $i+1$ over the days exhibiting the largest average hourly deviation between $\boldsymbol{\tilde{\lambda}}^i$ and $\boldsymbol{\hat{\lambda}}^i$,
rather than over randomly sampled days as in Algorithm~\ref{alg1},
adapting TSA to the observed optimality gap.

Critically, both Algorithms~\ref{alg1} and~\ref{alg2} exhibit the following key properties.
First, they are accompanied by a formal performance guarantee,
in the form of theoretically validated bounds on the objective function error incurred relative to the original full-scale GEP model~\eqref{full_model}.
Second, at each iteration, they yield a feasible solution for~\eqref{full_model},
i.e., the decision variable values employed to compute the objective function upper bound.
Finally, they offer two complementary and practical strategies to exploit the theoretical results of \cite{10037240} and \cite{klatzer2025towards}, thereby guiding the construction of the aggregated model toward exact temporal aggregation based on marginal cost estimates.

\section{Numerical Results}
\label{sec:results}
This section presents the numerical results.
Subsection~\ref{subsec:res_case_study} describes the case study, Subsection~\ref{subsec:res_marginal_costs} illustrates how marginal-cost-based clustering enables exact TSA, and Subsection~\ref{subsec:res_algorithms} compares the proposed solution algorithms.

\subsection{Case Study Description}
\label{subsec:res_case_study}

To ensure transparency and reproducibility, we generate the problem parameters from prescribed distributions as follows.

For each storage unit,
the charging and discharging efficiencies are set to $0.9$ and $1/0.9$, respectively,
while the associated operational costs are drawn from a uniform distribution over $[5,15]$~\euro/MWh. 
Capacity investment costs are sampled uniformly from $[4.5 \times 10^5, 5.5 \times 10^5]$~\euro/MW,
and storage capacity investments are either zero or constrained to $[0.25, 1]$~MW.
The stored energy is initialized to $0$~MWh, with the energy-to-power ratio fixed at $2$~h.
Among the generators, $20$\% are thermal units, $40$\% wind, and $40$\% solar.
Capacity investments in thermal and vRES units are either zero or constrained to $[0.5, 1]$~MW and $[0.2, 1]$~MW, respectively.
The investment costs are sampled uniformly from $[3 \times 10^6, 4 \times 10^6]$~\euro/MW for thermal units
and $[5 \times 10^5, 6 \times 10^5]$~\euro/MW for vRES.
The operational costs are set to $50$~\euro/MWh for thermal units and $1$~\euro/MWh for vRES, 
while the cost of non-supplied energy is $10^5$~\euro/MWh.
The GEP horizon is one year with hourly resolution ($|\boldsymbol{T}| = 8760$).
In the proposed algorithms, we set $I = 25$, $\zeta = 10$ and $\epsilon^{\mathrm{thr}} = 0.01\%$.

The capacity factor of thermal units is fixed at $1$,
whereas vRES capacity factors and the energy demand time series are derived by scaling real-world observations from the Austrian power grid in 2024,
available on the ENTSO-E Transparency Platform \cite{hirth2018entso}, and illustrated in Fig.~\ref{fig:input_data}.
The baseline demand profile shown in Fig.~\ref{fig:input_data} is scaled by a factor $0.05\left(|\boldsymbol{G}| + |\boldsymbol{N}|\right)$ in order to reflect the aggregate system size.
To generate heterogeneous production profiles across vRES units, uniform multiplicative noise in the range $[0.85, 1.15]$ is applied to the baseline wind and solar capacity factors depicted in Fig.~\ref{fig:input_data}.

The optimization models are implemented on an Intel i7 processor with 32~GB of RAM, using Gurobi~12.0.1.

\begin{figure}[t]
\centerline{\includegraphics[scale=0.6]{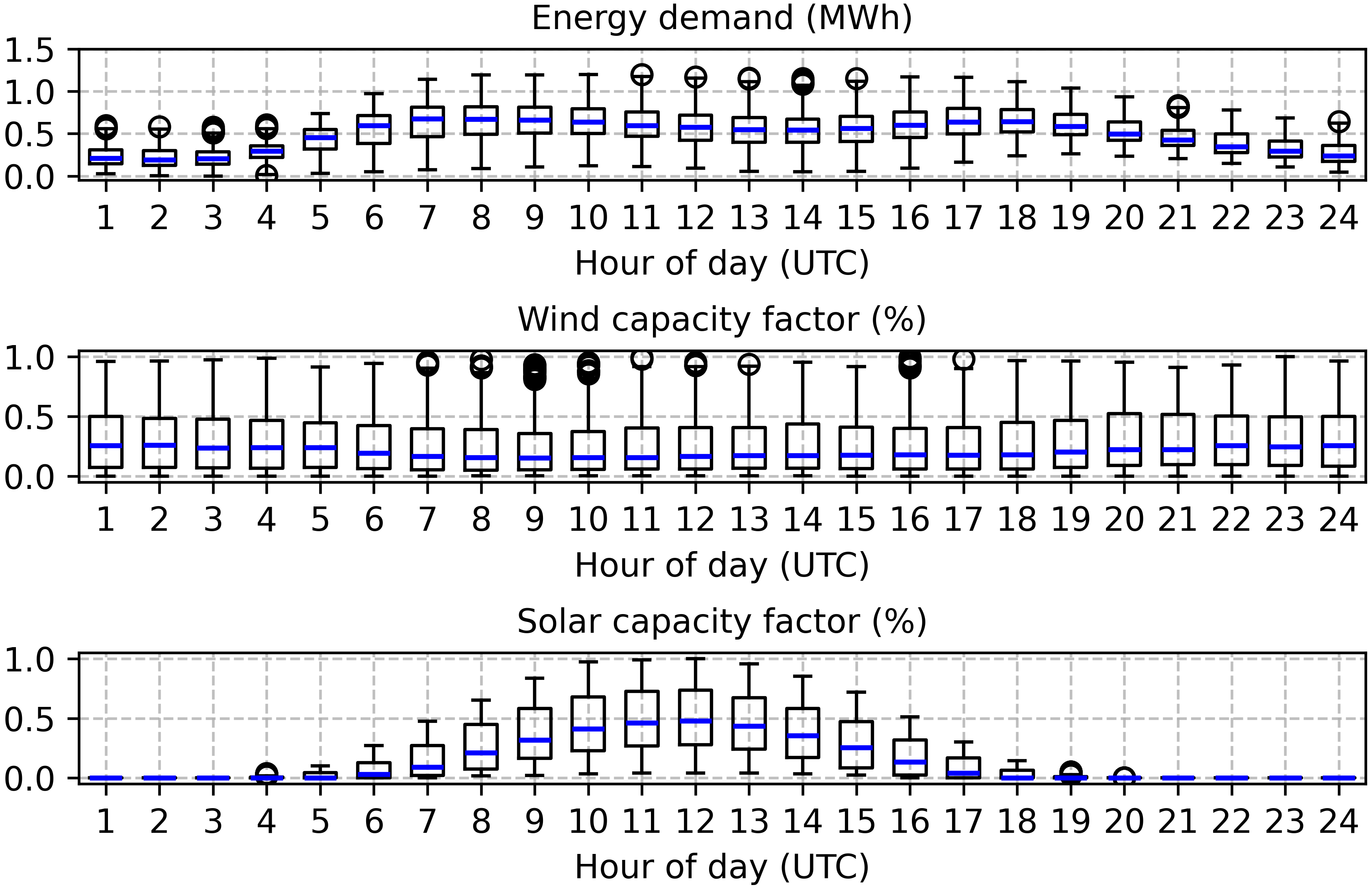}}
\caption{Boxplots of the hourly input data. Each boxplot characterizes the distribution of values observed at a given hour: the box spans the interquartile range, the blue line indicates the median, the whiskers extend to the 10th and 90th percentiles, and the outliers are shown as individual points.}
\label{fig:input_data}
\end{figure}

\subsection{Using Marginal Costs for Exact Time Series Aggregation}
\label{subsec:res_marginal_costs}

\begin{figure}[t]
\centerline{\includegraphics[scale=0.61]{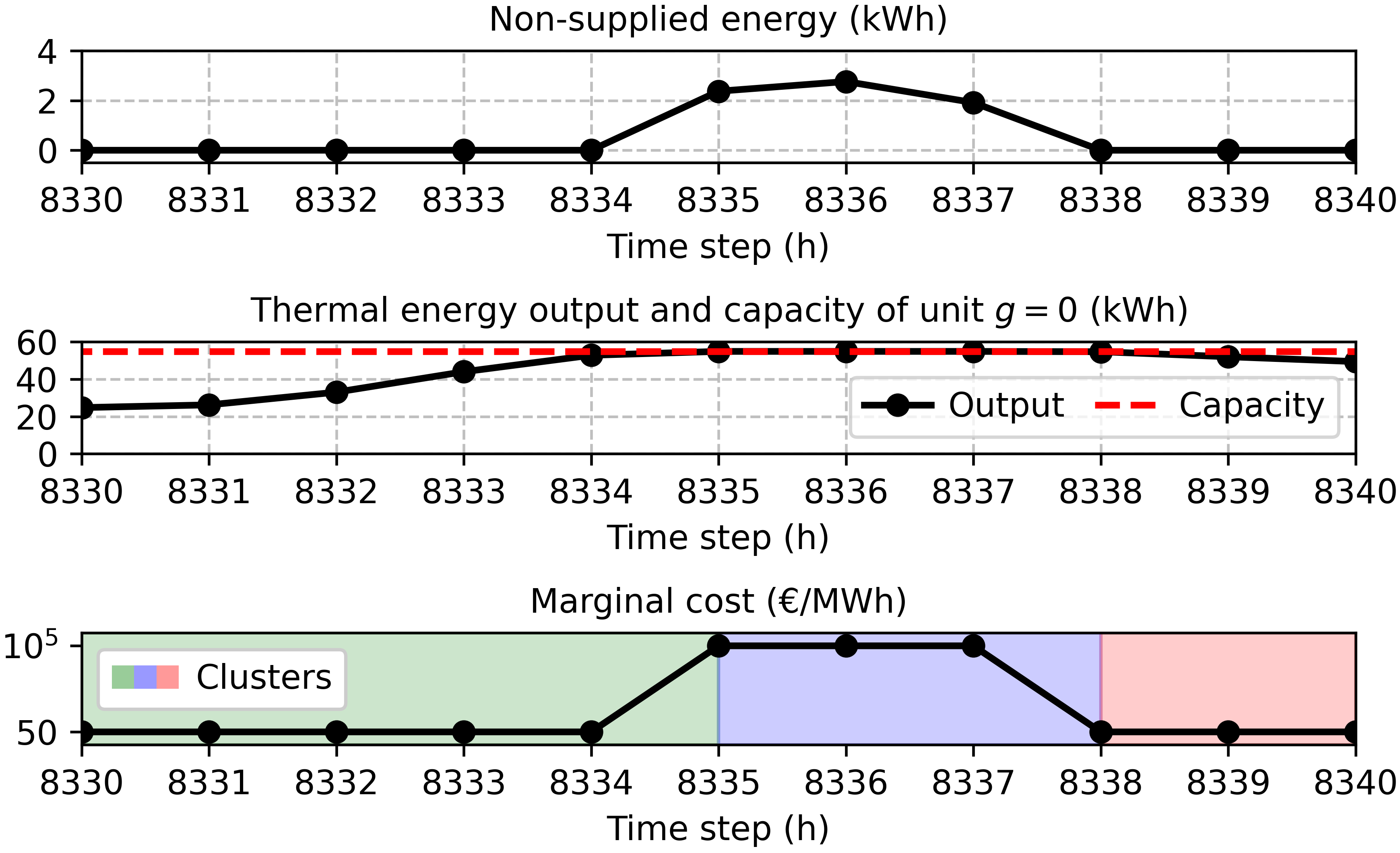}}
\caption{Example of marginal-cost-based time series aggregation.}
\label{fig:marginal_cost_clustering_example}
\end{figure}

Fig.~\ref{fig:marginal_cost_clustering_example} illustrates how marginal costs can be used to identify the active constraints in the LP relaxation of the full-scale GEP model~\eqref{full_model}.
Specifically, the results shown in the figure are obtained by setting $|\boldsymbol{G}| = 5$ and $|\boldsymbol{N}| = 0$,
thereby limiting the available sources of flexibility in the considered power system.
This configuration is deliberately chosen in this example to induce the activation of the power generation constraints \eqref{full_model:gen_lim} in the full-scale GEP model~\eqref{full_model}, resulting in the occurrence of non-supplied energy.

As shown in the figure, when the marginal generator (the only available thermal unit) reaches its maximum capacity (at $t=8335$), the marginal cost increases from $50\,\text{\euro}$ (its operating cost) to $10^5\,\text{\euro}$, corresponding to the non-supplied energy penalty. 
Conversely, as demand decreases and the generator operates strictly within its feasible region (starting from $t=8338$), the marginal cost decreases from $10^5\,\text{\euro}$ back to $50\,\text{\euro}$.
Hence, marginal costs alone constitute a sufficient statistic to identify the active constraint sets of the GEP problem,
specifically whether the upper limit on thermal power generation or the lower limit on non-supplied energy is binding.
This example illustrates that using marginal costs as clustering features for TSA yields an aggregated model that preserves the full-scale active constraints, thereby resulting in exact temporal aggregation, as demonstrated in \cite{10037240}.
In this sense, the informational value embedded in the system’s marginal costs is effectively unlocked, enabling exact TSA for the GEP problem under study.

\begin{figure}[t]
\centerline{\includegraphics[scale=0.51]{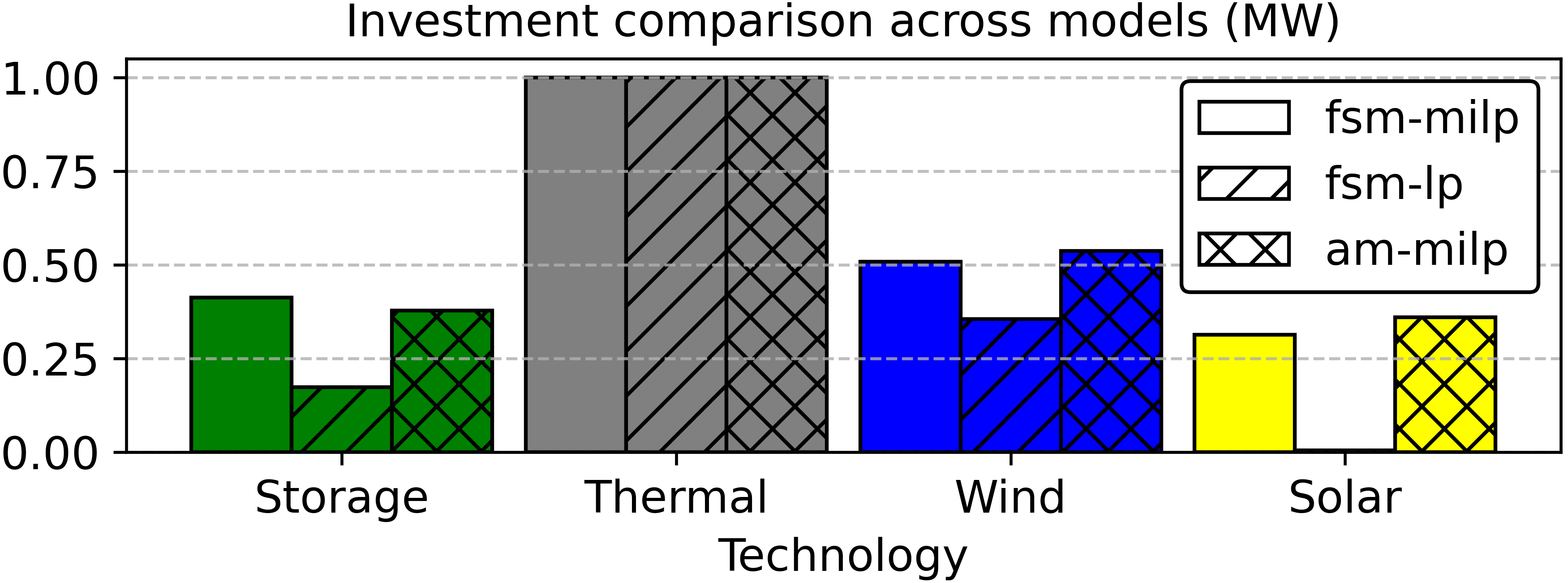}}
\caption{Comparison of investment decisions obtained with the full-scale MILP model (fsm-milp), its LP relaxation (fsm-lp), and the proposed aggregated MILP model (am-milp).}
\label{fig:inv_comparison_models}
\end{figure}

When addressing the MILP GEP model~\eqref{full_model}, marginal costs obtained from its LP relaxation are used as proxies for the true marginal costs within the proposed solution algorithms (see Subsection~\ref{subsec:algorithms}). 
Clearly, the LP relaxation could, in principle, serve directly as an approximation of the original MILP GEP model rather than solely for estimating marginal costs.
However, Fig.~\ref{fig:inv_comparison_models} shows that, for an illustrative stylized system comprising both storage ($|\boldsymbol{N}| = 1$) and generation ($|\boldsymbol{G}| = 10$) units, the combination of TSA based on marginal cost estimates from the LP relaxation with the construction of an aggregated model that preserves the nonconvex investment constraints, as proposed, provides a substantially more accurate approximation of the MILP GEP model’s optimal investment decisions than the conventional LP relaxation.

\subsection{Comparative Performance of the Proposed Algorithms}
\label{subsec:res_algorithms}

\begin{figure}[t]
\centerline{\includegraphics[scale=0.61]{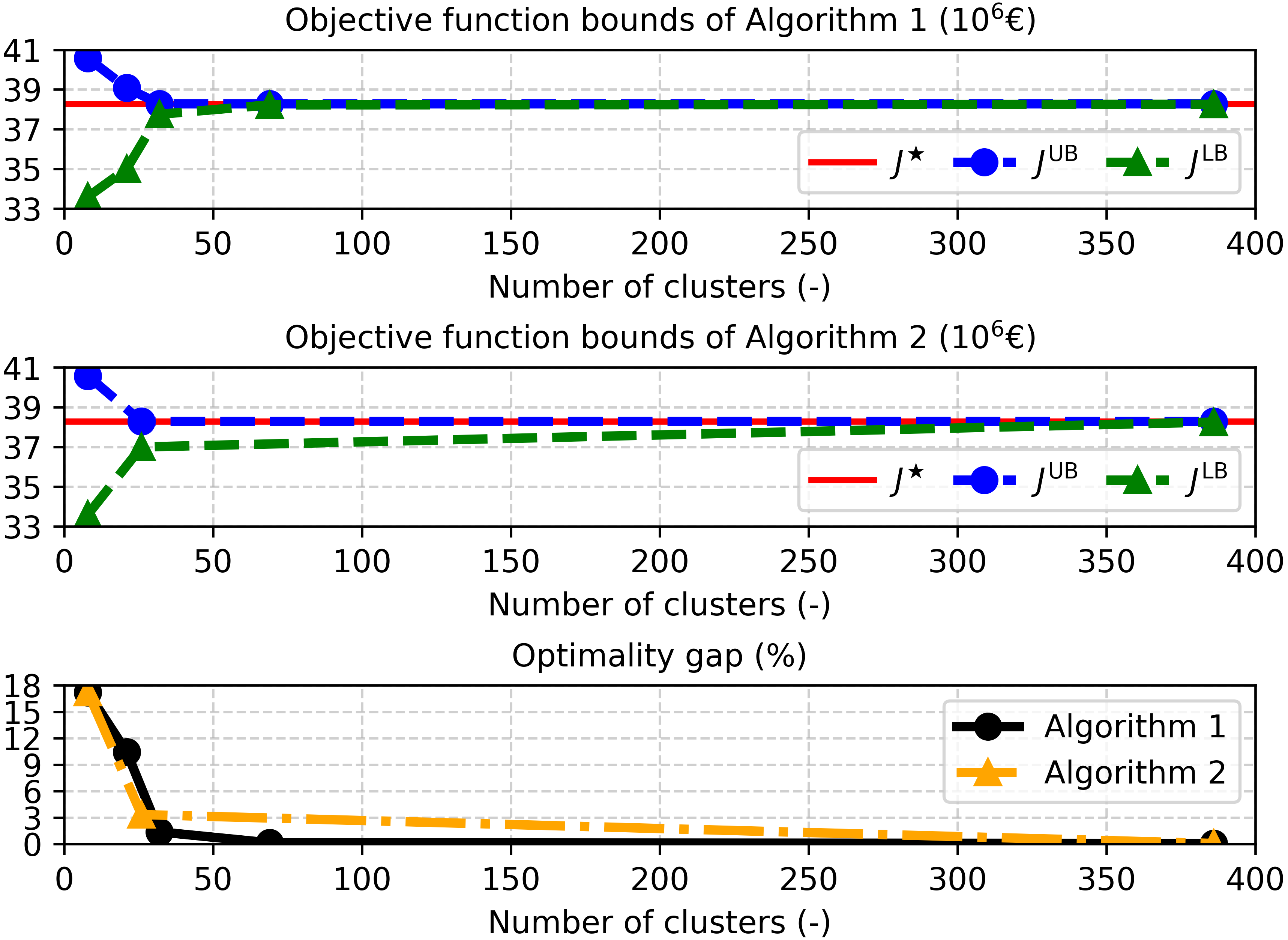}}
\caption{Performance comparison of the proposed Algorithms \ref{alg1} and \ref{alg2}.}
\label{fig:bounds}
\end{figure}

To illustrate the performance of the proposed solution algorithms, Fig.~\ref{fig:bounds} reports the convergence of their objective function bounds for a GEP model with $|\boldsymbol{G}| = 100$ and $|\boldsymbol{N}| = 10$.
Specifically, in the first two subplots of the figure, the lower and upper bounds on the optimal objective function value of the full-scale GEP model~\eqref{full_model}, obtained at each iteration of the proposed solution algorithms, are represented by triangles and circles, respectively.

As detailed in Subsection~\ref{subsec:algorithms}, both algorithms execute the same TSA and optimization steps during the first iteration.
From the second iteration onward, a key distinction emerges:
Algorithm~\ref{alg2} leverages the upper bound information obtained in the previous iteration to refine marginal cost estimates, whereas Algorithm~\ref{alg1} does not.
This difference is evident in Fig.~\ref{fig:bounds}, where Algorithm~\ref{alg2} achieves a $\sim3\%$ optimality gap by the second iteration, while Algorithm~\ref{alg1} exhibits a $\sim10\%$ gap.
Notably, both algorithms ultimately converge to over $99.99\%$ accuracy using $386$ clusters, corresponding to a reduction of more than $95\%$ in the temporal dimension of the GEP model.
However, Algorithm~\ref{alg1} requires five iterations to reach convergence,
whereas the adaptive Algorithm~\ref{alg2} converges in only three iterations, achieving the same dimensionality reduction with lower computational effort.

These results indicate that the computational advantage observed for Algorithm~\ref{alg2} over Algorithm~\ref{alg1} does not arise from achieving a more accurate temporal aggregation for a given number of clusters. Instead, it results from exploiting the information embedded in the upper bound computation, as described in Subsection~\ref{subsec:algorithms}, to accelerate the identification of the minimum number of clusters required to achieve the prescribed optimality gap.

\begin{table}[t]
\renewcommand{\arraystretch}{1.2}
\caption{Comparison between the runtimes of the full-scale model, Algorithm \ref{alg1}, and Algorithm \ref{alg2} as $|\boldsymbol{N}|$ and $|\boldsymbol{G}|$ increase. Relative differences with respect to the full-scale model are reported in brackets. The symbol $\infty$ indicates intractable instances of the full-scale model.}
\label{tab:comp_res}
\centering
\resizebox{8.8cm}{!}{\begin{tabular}{|c|c|c|c|c|}
\hline
\multirow{2}{*}{$|\boldsymbol{N}|$} & \multirow{2}{*}{$|\boldsymbol{G}|$} & \multicolumn{3}{c|}{\textbf{Runtime} (min)}\\
\cline{3-5}
& & \textbf{Full-scale model} & \textbf{Algorithm \ref{alg1}} & \textbf{Algorithm \ref{alg2}} \\
\hline
$50$ & $50$ & $2.1$ & $2.9$ ($\boldsymbol{+38\%}$) & $2.4$ ($\boldsymbol{+14\%}$)\\
$250$ & $250$ & $51.9$ & $22.4$ ($\boldsymbol{-57\%}$) & $19.6$ ($\boldsymbol{-62\%}$)\\
$750$ & $750$ & $432.8$ & $211.1$ ($\boldsymbol{-51\%}$) & $102.3$ ($\boldsymbol{-76\%}$)\\
$1500$ & $1500$ & $\infty$ & $291.9$ ($\boldsymbol{-\infty}$) & $144.2$ ($\boldsymbol{-\infty}$) \\
\hline
\end{tabular}}
\end{table}

Table~\ref{tab:comp_res} reports the computational performance of the proposed algorithms in comparison with conventional full-scale optimization.
As anticipated, Algorithm~\ref{alg2} consistently exhibits faster runtimes than Algorithm~\ref{alg1}. 
While for small instances of the GEP problem the TSA-based algorithms may incur higher computational costs relative to full-scale optimization,
their efficiency increases markedly as the number of generation and storage units,
and consequently the number of binary variables, grows.
Notably, both algorithms maintain tractability, whereas conventional full-scale optimization proves intractable, i.e., it fails to solve within the prescribed time of 24 hours.

\begin{table}[t]
\renewcommand{\arraystretch}{1.2}
\caption{
Percentage of total runtime devoted to the lower bound (LB) and upper bound (UB) computations in Algorithms \ref{alg1} and \ref{alg2}.
The LB computation includes TSA and the solution of the temporally aggregated model,
while the UB computation consists of solving the full-scale model with fixed investment decisions, as detailed in Section \ref{sec:methodology}.
}
\label{tab:comp_res_per_alg_step}
\centering
\begin{tabular}{|c|c|c|c|c|c|}
\hline
\multirow{3}{*}{$|\boldsymbol{N}|$} & \multirow{3}{*}{$|\boldsymbol{G}|$} & \multicolumn{4}{c|}{\textbf{Percentage of total runtime}}\\
\cline{3-6}
& & \multicolumn{2}{c|}{\textbf{Algorithm \ref{alg1}}} & \multicolumn{2}{c|}{\textbf{Algorithm \ref{alg2}}} \\
\cline{3-6}
& & \textbf{LB} & \textbf{UB} & \textbf{LB} & \textbf{UB} \\
\hline
$50$ & $50$
& $19\%$ & $81\%$
& $19\%$ & $81\%$\\
$250$ & $250$
& $31\%$ & $69\%$
& $29\%$ & $71\%$\\
$750$ & $750$
& $33\%$ & $67\%$
& $32\%$ & $68\%$\\
$1500$ & $1500$
& $38\%$ & $62\%$
& $34\%$ & $66\%$\\
\hline
\end{tabular}
\end{table}

Finally, Table \ref{tab:comp_res_per_alg_step} reports the percentage of the total runtime devoted to computing the objective function lower and upper bounds in the proposed solution algorithms.
The results indicate that, for small instances of the GEP problem, the upper bound computation accounts for the majority of the total runtime in both algorithms, representing more than $80\%$ of the overall computational effort.
This behavior is expected, since the upper bound computation requires solving the full-scale GEP model with investment decisions fixed to those obtained from the temporally aggregated model.
Consequently, despite being linear, the resulting optimization problem remains computationally demanding due to its large temporal dimensionality.

As the size of the GEP problem increases, the relative computational burden gradually shifts toward the lower bound computation.
Specifically, the share of total runtime devoted to this algorithmic step increases from $19\%$ for the smallest GEP problem to approximately $36$--$38\%$ for the largest problem considered.
This behavior reflects the different scalability characteristics of the two computational tasks:
while the upper bound computation requires solving a linear optimization model, the lower bound computation involves solving a temporally aggregated MILP model whose computational complexity grows more rapidly with the size of the GEP problem.

Notably, the growth in computational effort associated with the lower bound derivation is less pronounced for Algorithm \ref{alg2} than for Algorithm \ref{alg1}.
By leveraging the marginal cost information obtained during the objective function upper bound computation, Algorithm \ref{alg2} can accurately infer, after only a limited number of iterations, the number of representative time periods required to attain the prescribed optimality gap.
Consequently, several intermediate iterations that would otherwise be performed by Algorithm \ref{alg1} can be avoided, thereby reducing the number of temporally aggregated MILP models solved for refining the objective function lower bound during the execution of Algorithm \ref{alg2}.

\section{Conclusion and Future Work}
\label{sec:conclusion}
This paper addresses the GEP problem, formulated as a MILP model incorporating intertemporal storage constraints.
Being NP-hard, the problem’s computational complexity increases sharply with the planning horizon and the number of generators.
To address this challenge, we propose a novel a posteriori TSA method that leverages marginal cost estimates to ensure the aggregated model retains the active constraints of the full-scale model, thereby targeting exact temporal aggregation.
The TSA method is embedded within solution algorithms that iteratively refine theoretically validated bounds on the maximum error introduced by the temporal aggregation, while producing a feasible GEP solution at each iteration.

Numerical results show the potential of marginal-cost-based clustering for exact TSA and the efficiency of the proposed algorithms, which restore tractability where the full-scale model proves intractable.

Future work will focus on refining the marginal cost estimation procedure, incorporating additional nonconvex dynamics into the GEP formulation, and explicitly accounting for network constraints in larger-scale instances of the GEP problem.
In particular, extending the GEP model to include network constraints would require reformulating the proposed marginal-cost-based clustering technique of Subsection \ref{subsec:marginal_cost_clustering} into a clustering approach based on locational marginal prices.
In this context, the prediction of locational marginal prices would depend not only on the system-wide marginal generation cost, as considered in the present study, but also on network congestion and transmission losses.
To this end, machine learning-based prediction methods, recently initiated in our previous work \cite{rybka2026machine}, will be further investigated to assess whether they can preserve a favourable balance between prediction accuracy and computational efficiency.
Moreover, the marginal cost estimation procedure proposed in Subsection \ref{subsec:algorithms} could be further investigated by replacing the current random sampling of days within the GEP horizon with a feature-informed sampling strategy.
Specifically, selecting days based on relevant system characteristics, such as demand profiles and vRES capacity factors, may improve the quality of the marginal cost estimates and potentially enhance the convergence of the proposed solution algorithms.

\section{AI Usage Disclosure}
AI was utilized for spelling and grammar checking.

\section*{Acknowledgments}
Funded by the European Union (ERC, NetZero-Opt, 101116212). Views and opinions expressed are however those of the authors only and do not necessarily reflect those of the European Union or the European Research Council. Neither the European Union nor the granting authority can be held responsible for them.

\endgroup
\end{document}